\documentclass{article}

\usepackage{arxiv}

\usepackage[utf8]{inputenc} 
\usepackage[T1]{fontenc}    
\usepackage{hyperref}       
\usepackage{url}            
\usepackage{booktabs}       
\usepackage{amsfonts}       
\usepackage{nicefrac}       
\usepackage{microtype}      
\usepackage{cleveref}       
\usepackage{graphicx}
\usepackage[numbers]{natbib}
\usepackage{doi}
\usepackage{subcaption}
\newtheorem{theorem}{Theorem}
\newtheorem{Proposition}[theorem]{Proposition}%

\newtheorem{Example}{Example}%
\newtheorem{remark}{Remark}%

\title{Exact Solutions for Nonlinear Partial Differential Equations: A Fusion of Classical Methods and Innovative Approaches}

\author{Noureddine Mhadhbi \thanks{Department of Mathematics, College of Sciences and Arts, King Abdulaziz University. P.O. Box 344, Rabigh Campus, 21911, Saudi Arabia.
		\texttt{nalmhadhbi@kau.edu.sa} }
	\And
	Sameh Gana \thanks{
		Department of Basic Sciences. Deanship of Preparatory Year and Supporting Studies,
		Imam Abdulrahman Bin Faisal University. P.O. Box 1982, Dammam, 34212, Saudi Arabia.
		\texttt{sbgana@iau.edu.sa}}
	\And
	Mazen Fawaz Alsaeedi \thanks{Department of Mathematics, College of Sciences and Arts, King Abdulaziz University. P.O. Box 344, Rabigh Campus, 21911, Saudi Arabia.
		\texttt{malsaeedi0022@stu.kau.edu.sa}}}
\begin{document}
	\maketitle
	\begin{abstract}
		This article demonstrates how variation of parameters can be successfully implemented in combination with other classical techniques, such as the method of characteristics, to derive novel classes of solutions to nonlinear partial differential equations(NLPDES) by considering specific initial conditions. This innovative approach offers the advantage of generating exact solutions. The results underscore this method's potential to address additional NLPDE classes.
		
	\end{abstract}


	\section{Introduction}
	
	Nonlinear partial differential equations are prevalent in many physical problems, such as solid mechanics, fluid dynamics, acoustics, nonlinear optics, plasma physics, and quantum field theory. They also find applications in chemical and biological systems and formulate the fundamental laws of nature. Within this broad spectrum, a particularly intriguing class of nonlinear partial differential equations known as soliton equations gives rise to physically attractive solutions known as solitons. These solitons have significantly contributed to the field of applied sciences. For a comprehensive understanding of these phenomena, refer to ~\citep{Polyanin2},~\citep {Faddeev},~\citep {Gana1},~\citep {Gana3},~\citep{Myint},~\citep{Zachmanoglou} and ~\citep {Rhee} and the references therein, which offer detailed insights from both the theoretical and experimental perspectives. 
	
	Pursuing suitable analytical methods to solve nonlinear partial differential equations is a central focus. 
	Among the most widely adopted techniques are the variational iteration method ~\citep{Wu}, the inverse scattering method ~\citep{Ablowitz}, the truncated expansion method ~\citep{Kudryashov}, the extended tanh-function method ~\citep{Fan}, 
	Jacobi elliptic method ~\citep{Shikuo}, the Backlund transformations ~\citep{Hong}, F-expansion method ~\citep{D.Zhang,Zhang}, the sine-cosine function method ~\citep{Wazwaz}, the (G'/G)-expansion method ~\citep{Wang},
	and various extensions.
	
	One of the valuable tools for solving certain types of PDEs is the method of characteristics ~\citep{Polyanin2}, ~\citep{Myint}, ~\citep{Zachmanoglou}, ~\citep {Rhee} and ~\citep {Brian}. It involves transforming a PDE into a set of ordinary differential equations along characteristic curves. The characteristic curves represent the paths along which the solution of the PDE remains constant. The method of characteristics is a powerful technique for solving first-order partial differential equations (PDEs), including linear first-order PDEs such as the transport equation or the linear advection equation. 
	
	The well-known classical method usually refers to the variation of parameters ~\citep{Mahouton}, ~\citep{Polyanin1}, ~\citep{Jovan}, ~\citep{Olver} and ~\citep{Kevorkian}. The variation of parameters is primarily a technique used for linear differential equations, both ordinary and partial. It involves finding a particular solution to a non-homogeneous equation by introducing a new function to replace a constant in the homogeneous solution. 
	Solving nonlinear partial differential equations (NLPDEs) can pose considerably greater complexity and demand a problem-specific approach since nonlinear equations lack the superposition properties present in linear equations. The approach relies on the particular structure and characteristics of the NLPDE being addressed.
	
 The variation of parameters method has been successfully applied to certain NLPDES. We can refer to  ~\citep{Mahouton} and ~\citep{Jovan} as interesting studies. \\
	Common examples of second-order equations that can be converted into first-order forms include various types of nonlinear wave equations, heat equations, and specific conservation laws. The exact procedure for this reduction may vary depending on the specific equation and the desired format for further analysis.
	
	However, the method of characteristics and the variation of parameters are two distinct methods used in different contexts. While these two methods have distinct applications, this study demonstrates that combining the classical techniques derives new solutions for NLPDEs with specific initial conditions.
	
	As an extension of a previous study ~\citep{Mhadhbi}, we introduced new solutions to NLPDEs.
	In this study, we consider the classes of nonlinear partial differential equations of the form: 
	$$u_{tt}+a(x,t)u_{xt}+b(t)u_{t}=\alpha(x,t)+ G(u)(u_{t}+a(x,t)u_{x})e^{-\int b(t)dt}$$ and 
	$$u_{t}^{m}(u_{tt}+a(x,t)u_{xt})+b(t)u_{t}^{m+1}=e^{-(m+1)\int b(t)dt}(u_{t}+a(x,t)u_{x})
	F(u,u_{t}e^{\int b(t)dt}) $$
	
	Notably, some exceptional cases can arise. For example, we mention the NLPDEs recorded in ~\cite{Mahouton}, where the functions were restricted to one variable. 
	
	The remainder of this paper is organized as follows. In Section 2, we apply our methodology to the first class of reducible second-order partial differential equations to determine the exact solutions of nonlinear partial differential equations of the first type. 
	
	Section 3 delves into the second class of reducible nonlinear partial differentiable equations. Based on these results, a new class of solutions was derived. We demonstrate the application of the proposed method using concrete examples to demonstrate its viability and efficiency. Using Mathematica algorithms, relevant numerical representations were exhibited in each example to show the pertinence of obtained analytical solutions. Finally, Section 4 concludes the paper.

	\section{First class of reducible nonlinear partial differential equations}
	\subsection{Description of the method and construction of the general solutions}
	
	We consider the first class of nonlinear second-order partial differential equations compilable in the following general form:
	\begin{equation}\label{eq1}
		u_{tt}+a(x,t)u_{xt}+b(t)u_{t}=\alpha(x,t)+ G(u)(u_{t}+a(x,t)u_{x})e^{-\int b(t)dt}
	\end{equation}
	where $u$ denotes a function of $(x,t)\in \mathbb{R}^{2}$.
	\\
	First, we solve the characteristic equation  $$\left\{
	\begin{array}{ccc}
		\frac{d}{dt}(x(t))&=&a(x(t),t) \\
		x(0)&=&x_{0}.
	\end{array}
	\right. $$ 
	Then (\ref{eq1}) can be rewritten as \begin{equation}u_{tt}+a(x(t),t)u_{xt}+b(t)u_{t}=\alpha(x(t),t)+ G(u)(u_{t}+a(x(t),t)u_{x})e^{-\int b(t)dt}.
	\end{equation}
	Multiplying both sides of (2) by $e^{\int b(t)dt}$, we get 
	\begin{equation}\label{eq2}\frac{d}{dt}(u_{t}(x(t),t) e^{\int b(t)dt})=\alpha(x(t),t)e^{\int b(t)dt}+ G(u)(u_{t}+a(x(t),t)u_{x}).
	\end{equation}
	The nonlinear second-order partial differential equation (\ref{eq1}) can be solved easily if we assume that 	$$u_{t}(x,t)=(H(t)+K(u))e^{-\int b(t)dt}.$$ where $H$ and $K$ are differentiable functions of $t$ and $u$ respectively.
	
	Then, we differentiate to obtain 
	\begin{equation}\label{eq3}\frac{d}{dt}(u_{t}(x(t),t) e^{\int b(t)dt})=H'(t)+K'(u)(u_{t}+a(x(t),t)u_{x}).\end{equation}
	Substituting (\ref{eq3}) into (\ref{eq2}), we find that:
	\begin{equation}\label{eq4} H'(t)=\alpha(x(t),t)e^{\int b(t)dt}
	\end{equation} and 
	\begin{equation}\label{eq5} K'(u)=G(u).
	\end{equation}
	Therefore, we arrive at the following results:
	\begin{Proposition}
		The second order partial differential equation (\ref{eq1}) can be reduced to the first order differential equation  $$u_{t}(x,t)=(H(t)+K(u))e^{-\int b(t)dt},$$ where the functions $H$ and $K$ are the general solutions of 
		$H'(t)=\alpha(x(t),t)e^{\int b(t)dt}$ and 
		$K'(u)=G(u).$
	\end{Proposition}
	\begin{remark}
		Let $G=u^{n}$, where $n$ is a non zero positive
		integer.\end{remark} Then, the second order partial differential equation (\ref{eq1}) becomes
	$$u_{tt}+a(x,t)u_{xt}+b(t)u_{t}=\alpha(x,t)+ u^{n}(u_{t}+a(x,t)u_{x})e^{-\int b(t)dt}$$
	Applying equations (\ref{eq4}) and (\ref{eq5}), we get an Abel equation of the form
	$$u_{t}(x,t))=(H(t)+\frac{u^{n+1}}{n+1})e^{-\int b(t)dt}.$$
	A comprehensive compilation of integrable Abel equations can be found in ~\citep{D.Zwillinger,G.M.Murphy,Kamke} and ~\citep{Panayotounakos}.
	\subsection{Application}
	\begin{Example}  Let $a(x,t)=x$, $\alpha(x,t)=xe^{-t}$, $b(x,t)=1$ and $G(u)=2u$
		\begin{equation}\label{eq6} u_{tt}+xu_{xt}=xe^{-t}+2u(u_{t}+xu_{x})e^{-t}.  \end{equation}
		with the initial conditions $u(x,0)=1$ and $u_{t}(x,0)=x+1$.
	\end{Example}
	
	\textbf{Solution:}
	
	We solve the characteristic equation
	$$\left\{
	\begin{array}{ccc}
		\frac{dx(t)}{dt}&=&x\\
		x(0)&=&x_{0}
	\end{array}
	\right. $$ 
	which leads to $x(t)=x_{0}e^{t}.$\\
	The functions $H$ and $K$ are general solutions of
	\begin{center}
		\begin{tabular}{ccc}
			$H'(t)$&=& $x(t)e^{-t}e^{t}=x(t)$,\\
			$K'(u)$&=& $2u.$
		\end{tabular}
	\end{center}
	Then we get
	$$H(t)=x
	+C_{1},$$
	and $$K(u)=u^{2}+C_{2}.$$
	where $C_{1}$ and $C_{2}$ are arbitrary constants.
	
	The second-order partial differential equation (\ref{eq6}) is reduced to the first-order differential equation
	\begin{equation}\label{eq7} 
		u_{t}(x,t)=(u^{2}+x)e^{-t},
	\end{equation}
	with initial condition $u(x,0)=1$.\\
	The first-order differential equation (\ref{eq7}) is an Abel equation, which can be solved using various methods. For more details, refer to ~\citep{D.Zwillinger,G.M.Murphy,Kamke}.
	
\begin{figure}[h]
		\centering
		\begin{subfigure}[b]{0.5\textwidth}
			\centering
			\includegraphics[width=\textwidth]{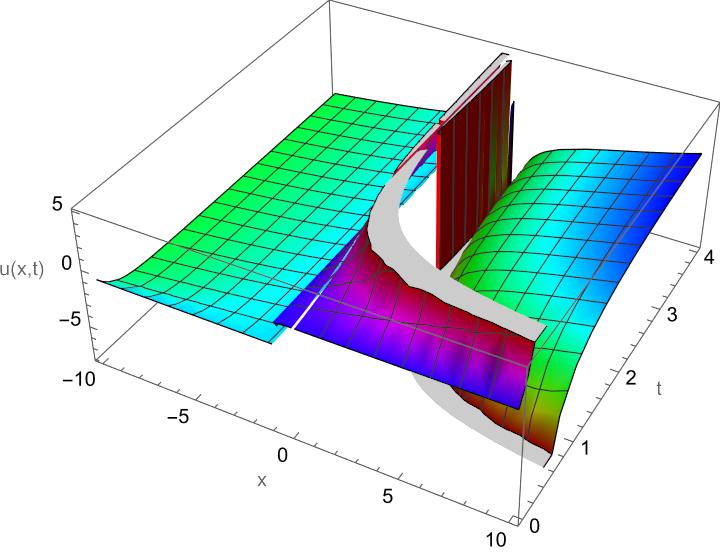}
			\caption{ }
		\end{subfigure}
		\hfill
		\begin{subfigure}[b]{0.4\textwidth}
			\centering
			\includegraphics[width=\textwidth]{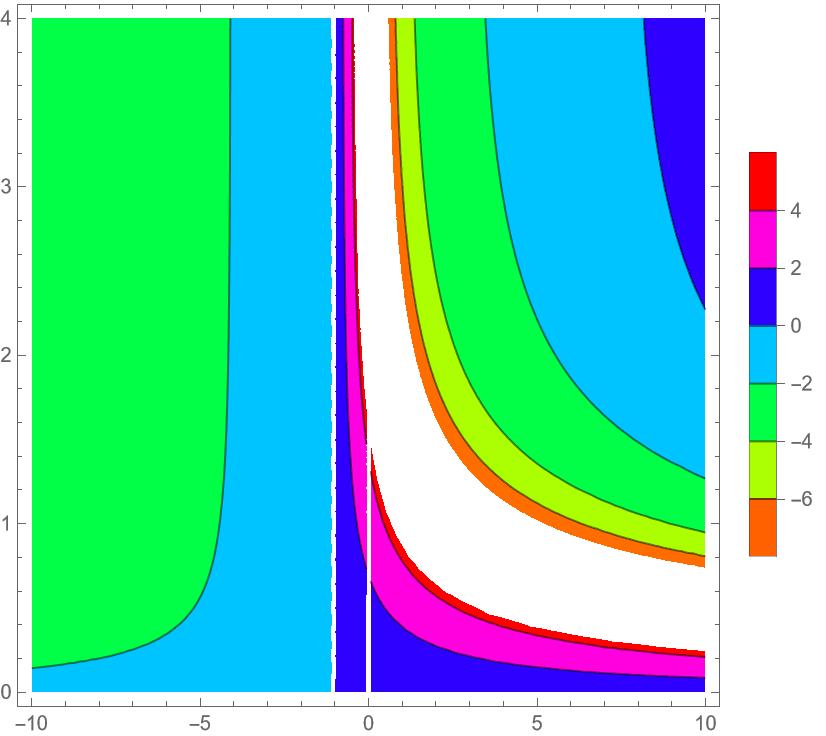}
			\caption{ }
		\end{subfigure}
		\hfill
		\begin{subfigure}[b]{0.5\textwidth}
			\centering
			\includegraphics[width=\textwidth]{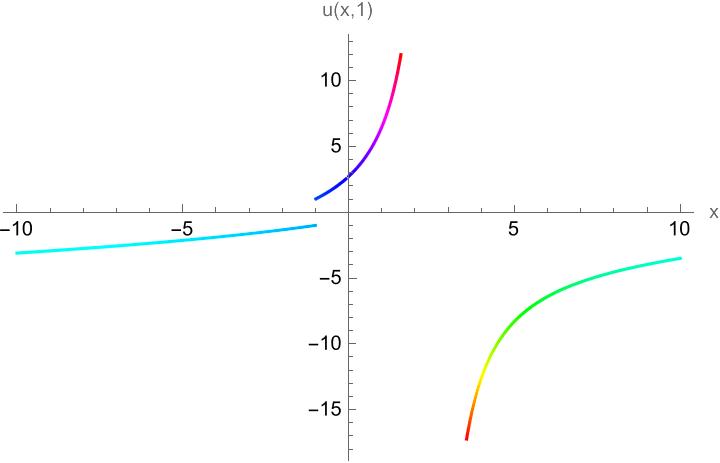}
			\caption{}
		\end{subfigure}
		\caption{The profile of the solutions in  (\ref{eq6}) with $u(x,0)=x$ and $u_{t}(x,0)=x+1$: (a) and (b) 3D and Contour plots  with $-10\leq x\leq10$ and $0\leq t\leq 4$, (c) 2D plot at $t=1$.}
		\label{fig1}
	\end{figure}

	Using initial condition $u(x,0)=1$ and $u_{t}(x,0)=x+1$ , we obtain explicit solutions of (\ref{eq6})
	$$u(x,t)=\sqrt{x}\, \tan \left(\sqrt{x}\, {\mathrm e}^{-t} \left(-1+\frac{e^{t} \left(\sqrt{x}-arcos \left(-\frac{\sqrt{x}}{\sqrt{1+x}}\right)\right)}{\sqrt{x}}\right)\right).$$
	Visualizing the precise solutions obtained by Mathematica algorithms (figure \ref{fig1}).
	By plotting solution profiles at different values of $t$, we observe the characteristics of several solutions of (\ref{eq6}) with initial conditions $u(x,0)=1$ and $u_{t}(x,0)=x+1$. As a result, these solutions develop singularities at certain values of $x$ and $t$. Note that despite the smoothness of the initial data, the spontaneous singular behavior in the solutions must be due to the nonlinear term of the equation.\\
	Figure \ref{fig1} displays the 2D, 3D and contour plots of the solutions in (\ref{eq6}) within $-10\leq x\leq10$ and $0\leq t\leq 4$ for 3D and contour graphs, $t=1$ for 2D graph.
	
	\section{{Second Class of reducible nonlinear partial differential equations}}
		\subsection{Description of the method and construction of the general solutions}
	The second group of second-order partial differential equations is formulated as follows:
	\begin{equation}\label{eq8} u_{t}^{m}(u_{tt}+a(x,t)u_{xt})+b(t)u_{t}^{m+1}=e^{-(m+1)\int b(t)dt}(u_{t}+a(x,t)u_{x})
		F(u,u_{t}e^{\int b(t)dt})  \end{equation}
	As in the previous section, the second-order nonlinear partial differential equation (\ref{eq8}) can be readily solved if we suppose that 
	\begin{equation}\label{eq9}u_{t}(x,t)=K(u)e^{-\int b(t)dt}.\end{equation}
	where $u=u(x(t),t)$ and $K$ is a differentiable function.
	
	Then
	\begin{equation}\label{eq10}(u_{t}(x,t))^{m+1}=K^{m+1}(u)e^{-(m+1)\int b(t)dt}.\end{equation}
	Differentiate (\ref{eq10}) to obtain
	\begin{equation}\label{eq11}
		(m+1)u_{t}^{m}(u_{tt}+a(x,t)u_{xt})=(m+1)K^{m}(u)K'(u)(u_{t}+a(x,t)u_{x})e^{-(m+1)\int b(t)dt} \end{equation}
	$$ -b(t)(m+1)K^{m+1}(u)e^{-(m+1)\int b(t)dt}.$$
	Then 
	\begin{equation}\label{eq12}
		K^{m}(u)K'(u)(u_{t}+a(x,t)u_{x})e^{-(m+1)\int b(t)dt}-b(t)K^{m+1}(u)e^{-(m+1)\int b(t)dt }+b(t)u_{t}^{m+1}=\end{equation}
	$$e^{-(m+1)\int b(t)dt}(u_{t}+a(x,t)u_{x})F(u,u_{t}e^{\int b(t)dt}) .$$
	
	Substituting equations (\ref{eq9}), (\ref{eq10}) and (\ref{eq11}) into (\ref{eq12}), we find 
	$$K^{m}(u)K'(u)=F(u,K(u)).$$
	Therefore, the following statement holds
	\begin{Proposition}
		The second order nonlinear partial differential equation (\ref{eq8}) can be reduced to the first order differential equation  $$u_{t}(x,t)=K(u)e^{-\int b(t)dt}$$ 
		
		where the function $K$ is the general solution of 
		$K^{m}(u)K'(u)=F(u,K(u)).$ 
	\end{Proposition}
\subsection{Applications}
	\begin{Example} Let $m=0$, $b=0$ and $a=1$. Suppose that function $F$ satisfies $F(s,w)=w$.
		\begin{equation}\label{eq13} u_{tt}+u_{xt}=u_{t}^{2}+u_{x}u_{t}  \end{equation}
		
		with the initial conditions $u(x,0)=0$ and $u_{t}(x,0)=x^{2}$.
	\end{Example}
	
	\textbf{Solution:}
	
	Using the previous result, we find that the second-order nonlinear partial differential equation
	$$ u_{tt}+u_{xt}=u_{t}^{2}+u_{x}u_{t} $$
	can be reduced to the first-order differential equation 
	\begin{equation}\label{eq14}u_{t}(x,t)=K(u)
	\end{equation}  where $K$ is the general solution of 
	$K'(u)=K(u).$ 
	Taking $x(t)=t+x_{0}$, we obtain $K(u)=A(x_{0})e^{u}$ where $A$ is an arbitrary constant of integration.\\
	The differential equation (\ref{eq14}) takes the form $$u_{t}(x,t)=A(x_{0})e^{u}$$ and the exact solutions of 
	the second-order nonlinear partial differential equation (\ref{eq13}) are analytically determined and take the following form:
	$$u(x,t)=-\ln(F(x-t)+G(x))$$ where $F$ and $G$ are arbitrary functions.
	
	It follows from the initial conditions at $(x_{0},0)$  given by  $u(x,0)=0$ and $u_{t}(x,0)=x^{2}$ that the exact solutions of (\ref{eq13}) can be expressed explicitly as follows
	$$u(x,t)=-\ln(1-\frac{t^{3}}{3}+t^{2}x-tx^{2})$$
	
\begin{figure}[h]
		\centering
		\begin{subfigure}[b]{0.5\textwidth}
			\centering
			\includegraphics[width=\textwidth]{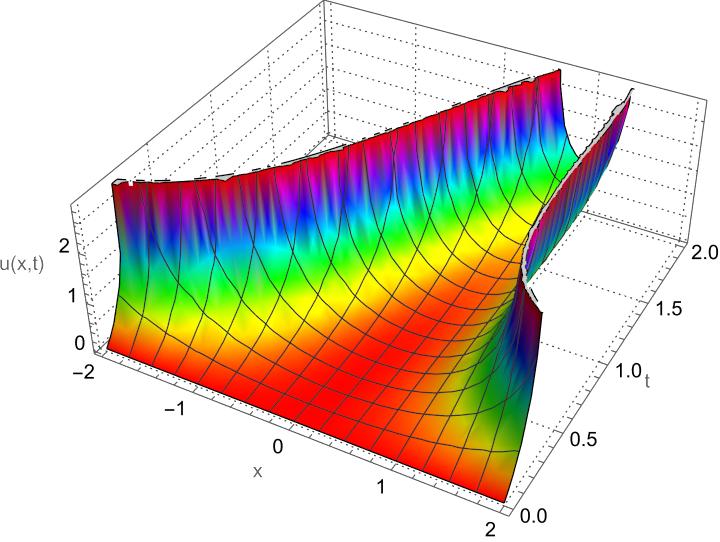}
			\caption{ }
		\end{subfigure}
		\hfill
		\begin{subfigure}[b]{0.4\textwidth}
			\centering
			\includegraphics[width=\textwidth]{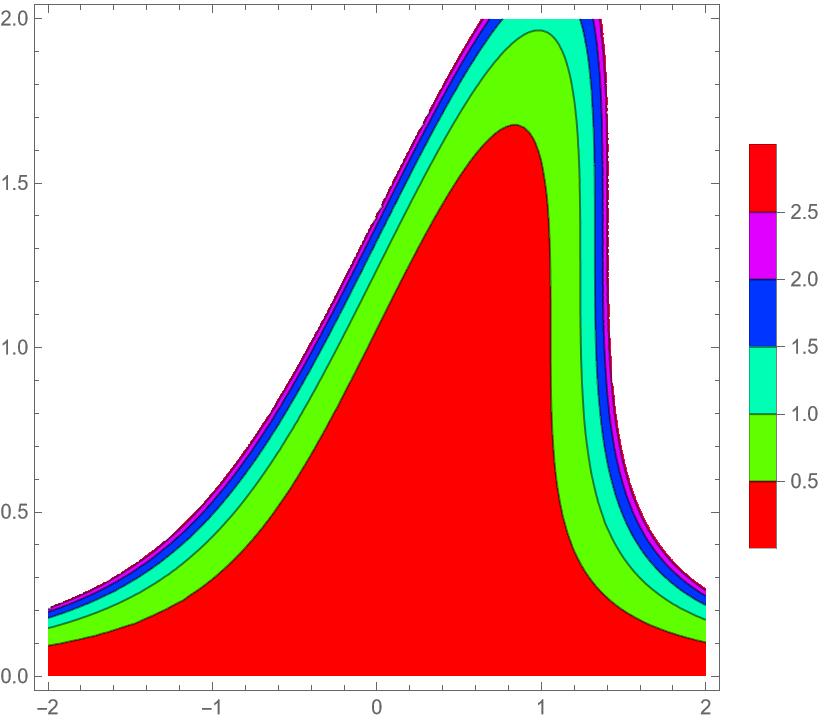}
			\caption{ }
		\end{subfigure}
		\hfill
		\begin{subfigure}[b]{0.5\textwidth}
			\centering
			\includegraphics[width=\textwidth]{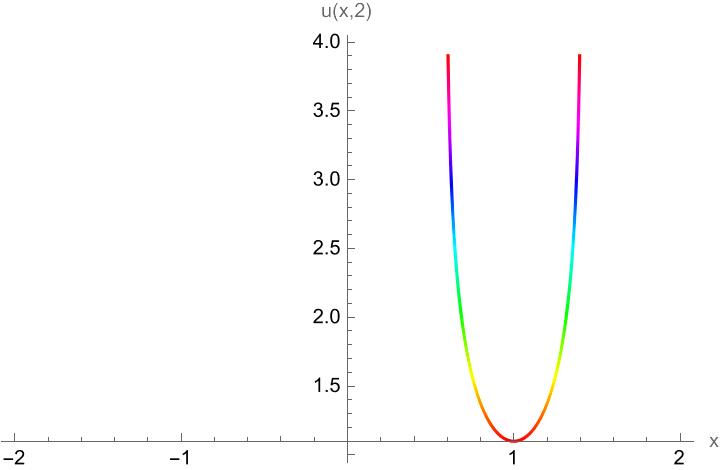}
			\caption{}
		\end{subfigure}
		\caption{The profile of the solutions in  (\ref{eq13}) with $u(x,0)=0$ and $u_{t}(x,0)=x^{2}$: (a) and (b) 3D and Contour plots  with $-2\leq x\leq2$ and $0\leq t\leq 2$, (c) 2D plot at $t=2$.}
		\label{fig2}
	\end{figure}

	Envisioning the precise solutions obtained by Mathematica (figure \ref{fig2}) and plotting solution profiles at different values of $t$, we have seen equations with smooth coefficients and initial data develop spontaneous singularities due to the nonlinearity of the equations. The solutions of (\ref{eq13}) break down at some values of $x$ and $t$, and no classical solution for the initial value problems exists beyond this point of breakdown. 
	
	Note that the nonlinear partial differential equation (\ref{eq13}) yields a more straightforward solution than the initial value problem in the previous example.
	
	Figure \ref{fig2} displays the 2D, 3D and contour plots of the solutions in (\ref{eq13}) within $-2\leq x\leq2$ and $0\leq t\leq 2$ for 3D and contour graphs, $t=2$ for 2D graph.
	
	\begin{Example} Let $m=0$, $b=1$ and $a=1$.
		
		Suppose the function $F$ satisfies $F(s,w)=w$.
		\begin{equation}\label{eq15} u_{tt}+u_{xt}+u_{t}=(u_{t}+u_{x})u_{t}  \end{equation} 
		
		with the initial conditions $u(x,0)=0$ and $u_{t}(x,0)=x^{2}$.
	\end{Example}
	\textbf{Solution:}
	
	Using the previous result, we find that the second-order nonlinear partial differential equation
	$$ u_{tt}+u_{xt}+u_{t}=(u_{t}+u_{x})u_{t}  $$
	can be reduced to the first-order differential equation 
	\begin{equation}\label{eq16}u_{t}(x,t)=e^{-t}K(u)
	\end{equation}  where $K$ is the general solution of 
	$K'(u)=K(u).$\\
	The differential equation (\ref{eq16}) takes the form $$u_{t}(x,t)=e^{-t}A(x_{0})e^{u}$$ and the exact solutions of 
	the second order nonlinear partial differential equation (\ref{eq15}) are analytically determined and take the following form:
	$$u(x,t)=-\ln(e^{-t}F(x-t)+G(x))$$ where $F$ and $G$ are arbitrary functions.
	
	It follows from the initial conditions at $(x_{0},0)$  given by  $u(x,0)=0$ and $u_{t}(x,0)=x^{2}$ that the exact solutions of (\ref{eq15}) can be expressed explicitly as follows
	$$u(x,t)=-\ln(e^{-t} t^2+e^{-t} x^2-2 e^{-t} t x-2 e^{-t} x+2 e^{-t} t+2 e^{-t}-x^2+2 x-1)$$
	
\begin{figure}[h]
		\centering
		\begin{subfigure}[b]{0.5\textwidth}
			\centering
			\includegraphics[width=\textwidth]{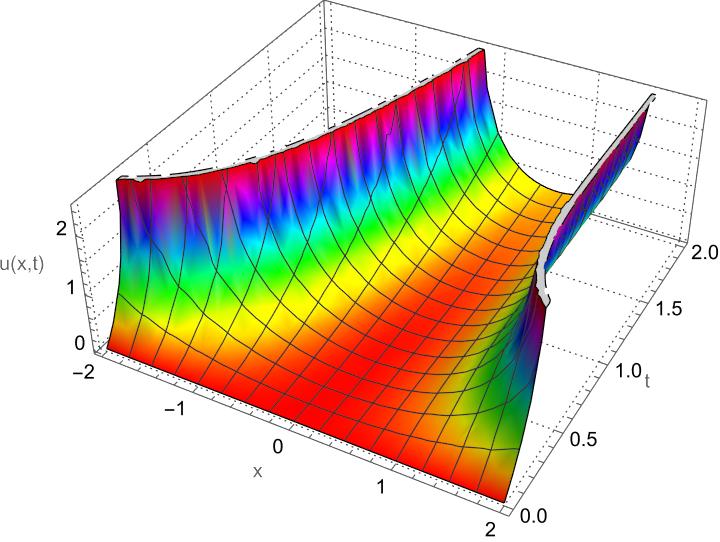}
			\caption{ }
		\end{subfigure}
		\hfill
		\begin{subfigure}[b]{0.4\textwidth}
			\centering
			\includegraphics[width=\textwidth]{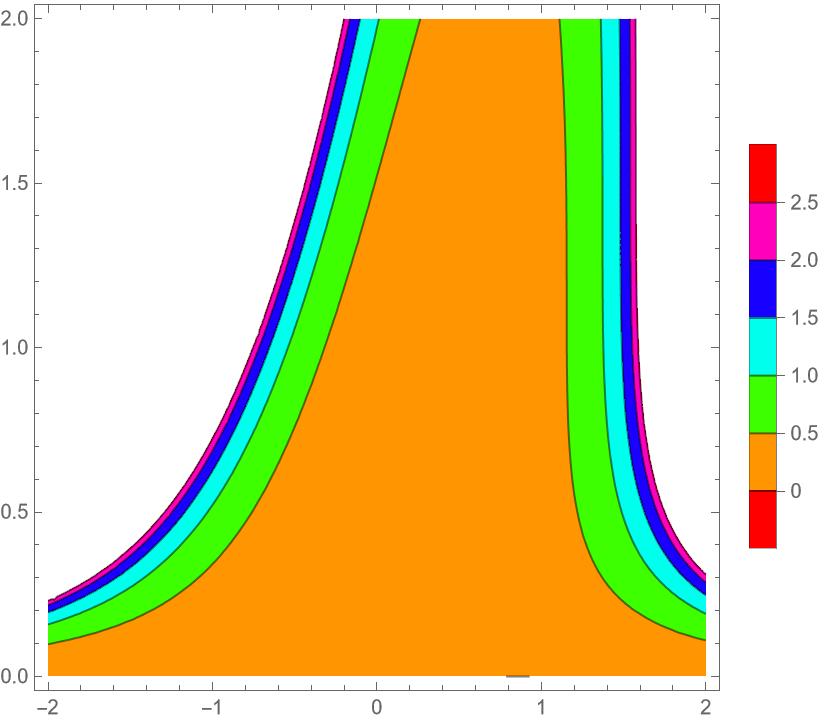}
			\caption{ }
		\end{subfigure}
		\hfill
		\begin{subfigure}[b]{0.5\textwidth}
			\centering
			\includegraphics[width=\textwidth]{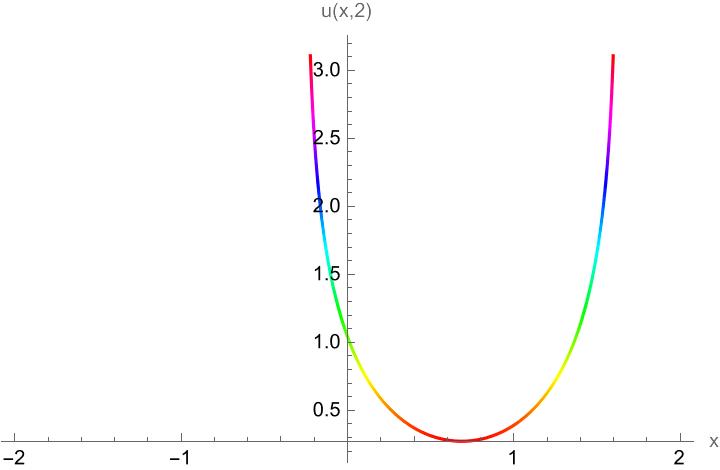}
			\caption{}
		\end{subfigure}
		\caption{The profile of the solutions in  (\ref{eq15}) with $u(x,0)=0$ and $u_{t}(x,0)=x^{2}$: (a) and (b) 3D and Contour plots  with $-2\leq x\leq2$ and $0\leq t\leq 2$, (c) 2D plot at $t=2$.}
		\label{fig3}
	\end{figure}

	When we envision the exact solutions of (\ref{eq15}) generated by Mathematica as depicted in Figure \ref{fig3} and create plots showing the solution profiles at various time points, we find that they deteriorate at specific values of both $x$ and $t$. Beyond this point, a classical solution is no longer viable for the initial value problems.
	
	Figure \ref{fig3} displays the 2D, 3D and contour plots of the solutions in (\ref{eq15}) within $-2\leq x\leq2$ and $0\leq t\leq 2$ for 3D and contour graphs, $t=2$ for 2D graph.

	\begin{Example} Let $m=0$, $b=0$ and $a=1$.
		Suppose that function $F$ satisfies $F(s,w)=s$.			
		\begin{equation}\label{eq17} u_{tt}+u_{xt}=(u_{t}+u_{x})u  \end{equation}
		
		with the initial conditions $u(x,0)=0$ and $u_{t}(x,0)=x^{2}$.\end{Example}
	
	\textbf{Solution:}	
	
	Using the previous result, we find that the second-order nonlinear partial differential equation
	$$u_{tt}+u_{xt}=(u_{t}+u_{x})u$$
	can be reduced to the first-order differential equation 
	\begin{equation}\label{eq18}u_{t}(x,t)=K(u)
	\end{equation}  where $K$ is the general solution of 
	$K'(u)=u.$\\
	The differential equation (\ref{eq18}) takes the form $$u_{t}(x,t)=\frac{1}{2}u^{2}+f(x-t)$$ where $f$ is an arbitrary function that leads to a Ricatti differential equation.
	
	Taking in account the initial conditions $u(x,0)=0$ and $u_{t}(x,0)=x^{2}$,  the differential equation (\ref{eq18}) becomes $$u_{t}(x,t)=\frac{1}{2}u^{2}+(x-t)^{2}$$
	The result was obtained using Mathematica code as a complicated function.  
	As in the previous examples, the nonlinearity of the partial differential equations produces singular behavior in the solutions.
	
\begin{figure}[h]
		\centering
		\begin{subfigure}[b]{0.5\textwidth}
			\centering
			\includegraphics[width=\textwidth]{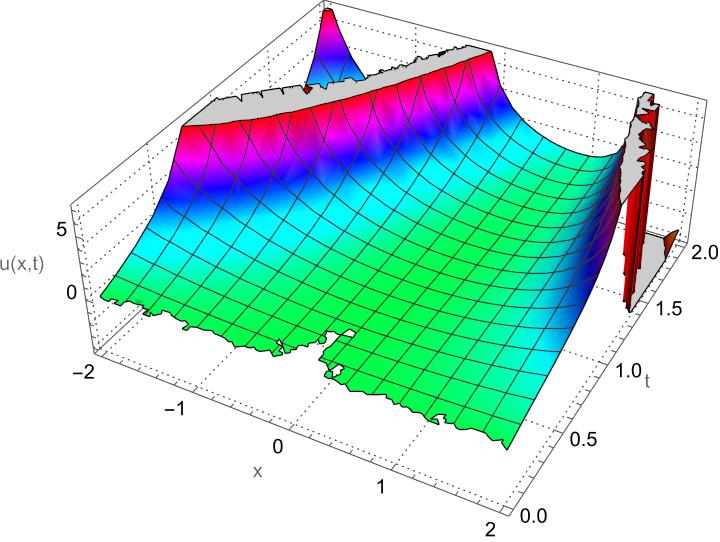}
			\caption{ }
		\end{subfigure}
		\hfill
		\begin{subfigure}[b]{0.4\textwidth}
			\centering
			\includegraphics[width=\textwidth]{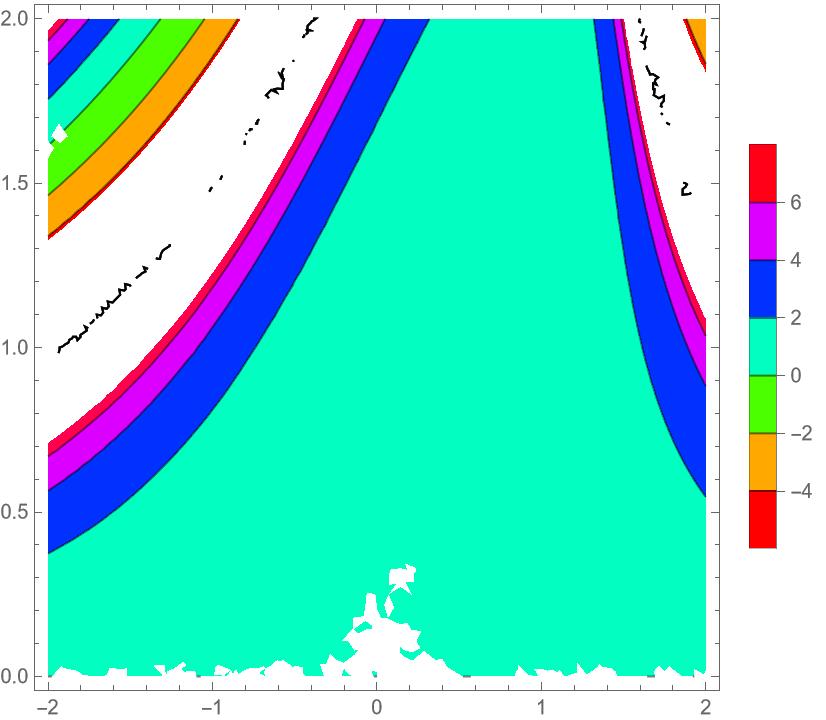}
			\caption{ }
		\end{subfigure}
		\hfill
		\begin{subfigure}[b]{0.5\textwidth}
			\centering
			\includegraphics[width=\textwidth]{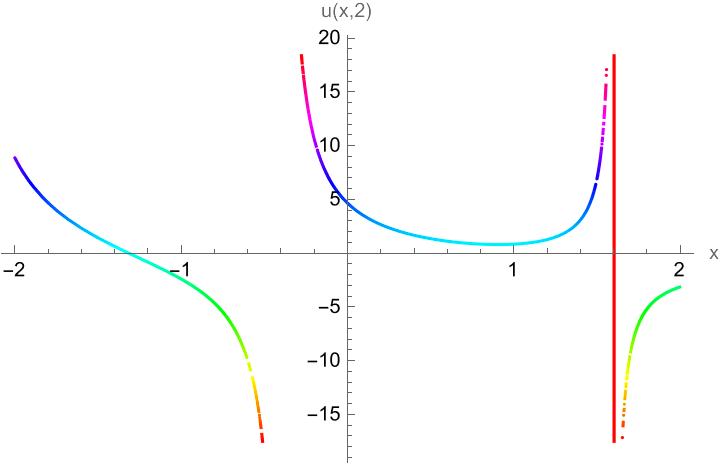}
			\caption{}
		\end{subfigure}
		\caption{The profile of the solutions in  (\ref{eq17}) with $u(x,0)=0$ and $u_{t}(x,0)=x^{2}$: (a) and (b) 3D and Contour plots  with $-2\leq x\leq2$ and $0\leq t\leq 2$, (c) 2D plot at $t=2$.}
		\label{fig4}
	\end{figure}

	Figure \ref{fig4} shows the 2D, 3D and contour plots of the solutions in (\ref{eq17}) within $-2\leq x\leq2$ and $0\leq t\leq 2$ for 3D and contour graphs, $t=2$ for 2D graph.
	
	\begin{Example} Let $m=0$, $b=0$ and $a=x$.
		
		Suppose that function $F$ satisfies $F(s,w)=s^{2}$.
		\begin{equation}\label{eq19} u_{tt}+xu_{xt}=(u_{t}+xu_{x})u^{2}  \end{equation} 	
		
		with the initial conditions $u(x,0)=x$ and $u_{t}(x,0)=\frac{x^{3}}{3}$.
		
	\end{Example}
	
	\textbf{Solution:}	
	
	Using the previous result, we find that the second-order nonlinear partial differential equation
	$$u_{tt}+xu_{xt}=(u_{t}+xu_{x})u^{2}$$
	can be reduced to the first-order differential equation 
	\begin{equation}\label{eq20}u_{t}(x,t)=K(u)
	\end{equation}  where $K$ is the general solution of 
	$K'(u)=u^{2}.$\\
	Differential equation (\ref{eq20}) takes the form of an Abel equation $$u_{t}(x,t)=\frac{1}{3}u^{3}+f(xe^{-t})$$ where $f$ is an arbitrary function.
	
	By applying the initial conditions $u(x,0)=x$ and $u_{t}(x,0)=\frac{x^{3}}{3}$, the exact solutions of (\ref{eq19})  are implicitly obtained by generating the Mathematica codes. 
	
\begin{figure}[h]
		\centering
		\begin{subfigure}[b]{0.5\textwidth}
			\centering
			\includegraphics[width=\textwidth]{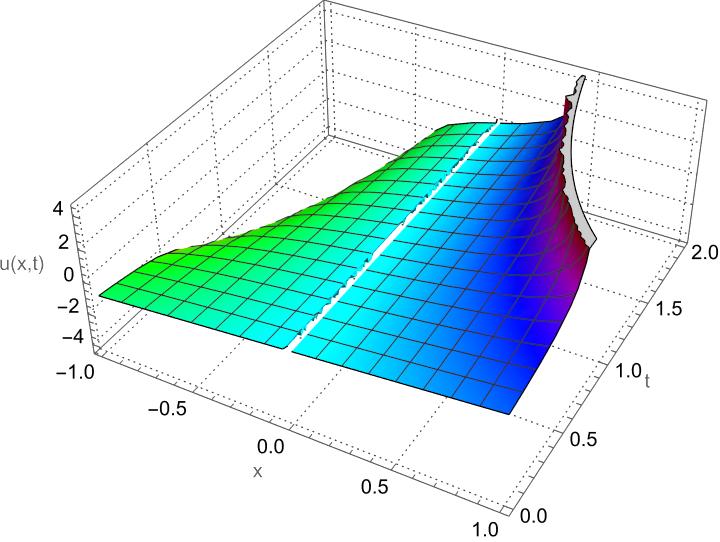}
			\caption{ }
		\end{subfigure}
		\hfill
		\begin{subfigure}[b]{0.4\textwidth}
			\centering
			\includegraphics[width=\textwidth]{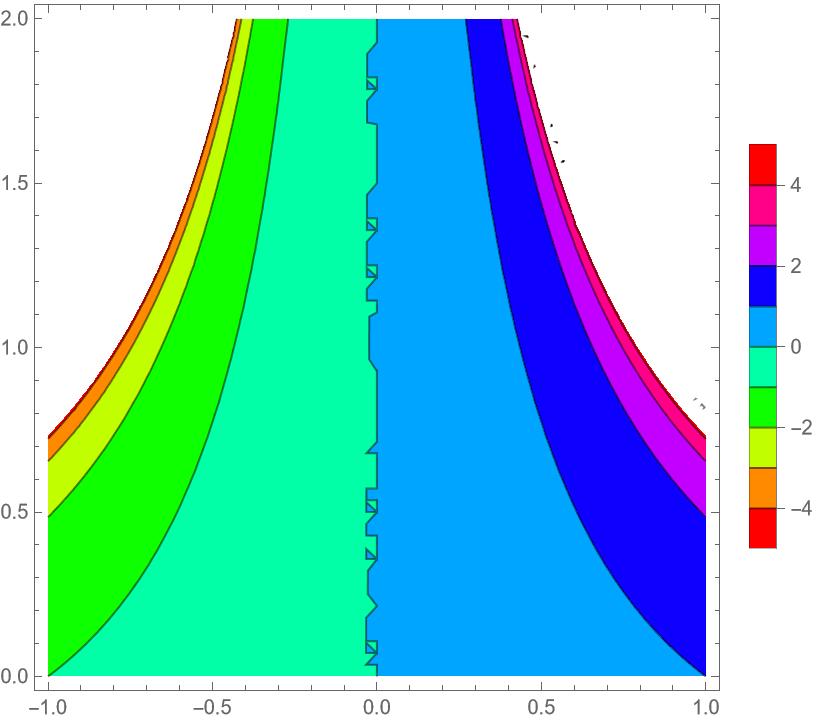}
			\caption{ }
		\end{subfigure}
		\hfill
		\begin{subfigure}[b]{0.5\textwidth}
			\centering
			\includegraphics[width=\textwidth]{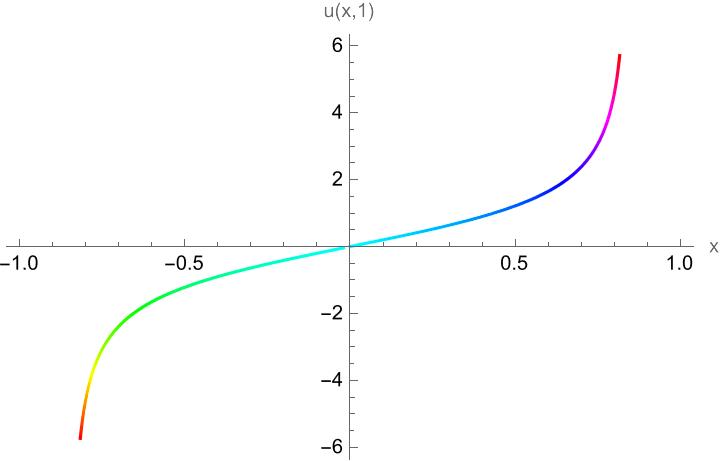}
			\caption{}
		\end{subfigure}
		\caption{The profile of the solutions in  (\ref{eq19}) with $u(x,0)=x$ and $u_{t}(x,0)=\frac{x^{3}}{3}$: (a) and (b) 3D and Contour plots  with $-1\leq x\leq1$ and $0\leq t\leq 2$, (c) 2D plot at $t=1$.}
		\label{fig5}
	\end{figure}

	Plotting the solution profiles for several values of $t$ (as depicted in Figure \ref{fig5}) shows that the solution breaks down at some points. 
	
	Figure \ref{fig5} shows the 2D, 3D and contour plots of the solutions in (\ref{eq19}) within $-1\leq x\leq1$ and $0\leq t\leq 2$ for 3D and contour graphs, $t=1$ for 2D graph.

	\begin{Example} Let $m=0$, $b(t)=\frac{1}{t}$.
		
		Suppose that function $F$ satisfies $F(s,w)=(\frac{w}{s})^{2} + 2\frac{w}{s}$.
		
		The second-order nonlinear partial differential equation (\ref{eq8}) becomes
		
		\begin{equation}\label{eq21}
			u_{tt}+a(x,t)u_{xt}+\frac{1}{t}u_{t}=\frac{1}{t}(u_{t}+a(x,t)u_{x})((\frac{tu_{t}}{u})^{2}+2\frac{tu_{t}}{u}) \end{equation} 
		
		with the initial conditions $u(x,1)=1$ and $u_{t}(x,1)=x^{2}$
	\end{Example}
	
	\begin{remark} Some interesting particular cases of (\ref{eq21}) can be formed. As an example, we mention the equation 
		\begin{equation}\label{eq22}u''(t)+\frac{1}{t}u'=tu'(\frac{u'}{u})^{2}+2\frac{u'^{2}}{u}\end{equation}  recorded in ~\citep{Mahouton} as equation (53).
		Equation (\ref{eq22}) is obtained from (\ref{eq21}) if we assume that $u$ is only a function of $t$.
		
	\end{remark}
	
	\textbf{Solution:}
	
	Using our previous result, we find that (\ref{eq21}) 
	can be reduced to the first-order differential equation 
	\begin{equation}\label{eq23}u_{t}(x,t))=K(u)\frac{1}{t}
	\end{equation}  where $K$ is the general solution of 
	$K'(u)=\frac{2}{u}K(u)+\frac{1}{u^{2}}K^{2}(u)$.\\
	We get $$K(u)=\frac{u^{2}}{A-u}.$$
	For $a=1$, (\ref{eq21}) takes the form 	
	\begin{equation}\label{eq24} u_{tt}+u_{xt}+\frac{1}{t}u_{t}=\frac{1}{t}(u_{t}+u_{x})((\frac{tu_{t}}{u})^{2}+2\frac{tu_{t}}{u})
	\end{equation} 
	Then the general solutions of (\ref{eq24}) are given by 
	$$u_{t}(x,t)=\frac{u^{2}}{f(x-t)-u}\frac{1}{t}$$ where
	$f$ denotes an arbitrary function.
	
	\begin{remark}
		Let $u(x,t)=t^{-1}B(x)$ be a family of solutions to equation (\ref{eq24}). \\
		If$f(x,t)=A$, we get 
		$$-\frac{A}{u(x,t)}-\ln(u(x,t))=\ln(t)+B(x)$$ which is an implicit solution of (\ref{eq24}) where $A$ an arbitrary constant and $B$ is a function of $x$.
	\end{remark}
	
	We checked the implicit solutions of (\ref{eq24}) by generating Mathematica codes and considering the initial conditions  $u(x,1)=1$ and $u_{t}(x,1)=x^{2}$.
	
\begin{figure}[h]
		\centering
		\begin{subfigure}[b]{0.5\textwidth}
			\centering
			\includegraphics[width=\textwidth]{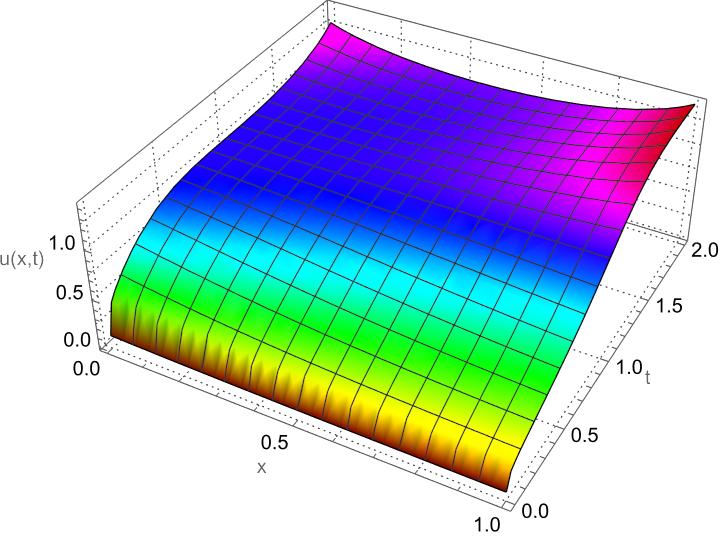}
			\caption{ }
		\end{subfigure}
		\hfill
		\begin{subfigure}[b]{0.4\textwidth}
			\centering
			\includegraphics[width=\textwidth]{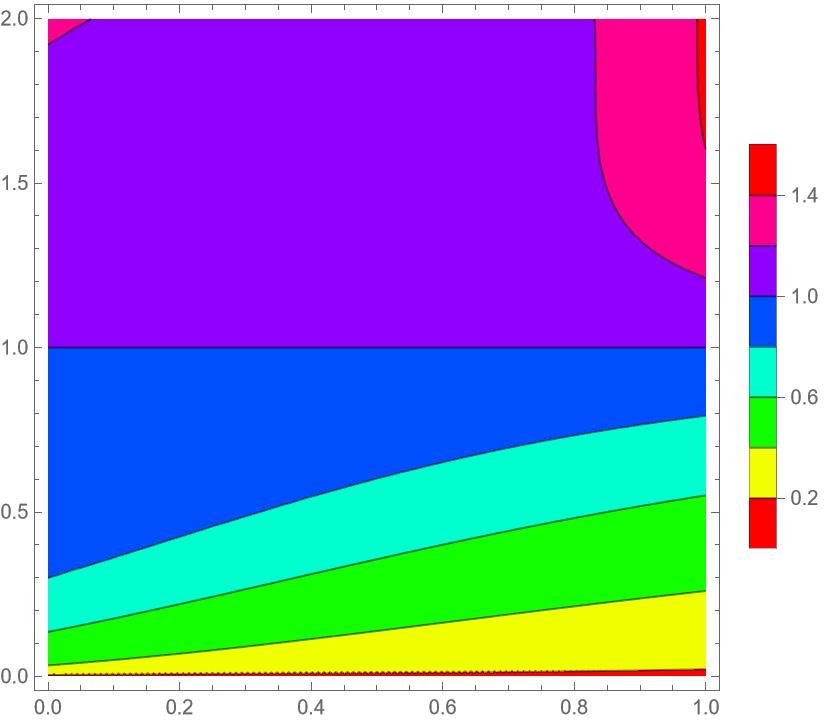}
			\caption{ }
		\end{subfigure}
		\hfill
		\begin{subfigure}[b]{0.5\textwidth}
			\centering
			\includegraphics[width=\textwidth]{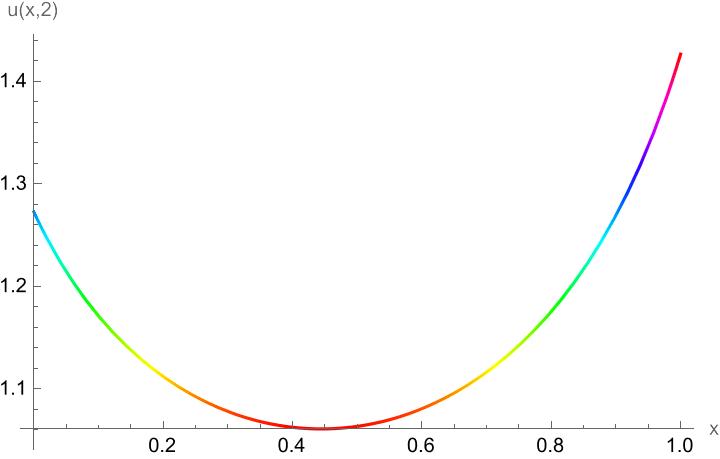}
			\caption{}
		\end{subfigure}
		\caption{The profile of the solutions in  (\ref{eq24}) with $u(x,1)=1$ and $u_{t}(x,1)=x^{2}$: (a) and (b) 3D and Contour plots  with $0\leq x\leq1$ and $0\leq t\leq 2$, (c) 2D plot at $t=2$.}
		\label{fig6}
	\end{figure}

	Figure \ref{fig6} shows the 2D, 3D and contour plots of the solutions in (\ref{eq24}) within $0\leq x\leq1$ and $0\leq t\leq 2$ for 3D and contour graphs, $t=2$ for 2D graph. 
	
	\begin{Example} Let $m=0$, $b(t)=\frac{2}{t}$ and $a=1$.
		Suppose that function $F$ satisfies $F(s,w)= w+w^{3}$.
		\begin{equation}\label{eq25} u_{tt}+u_{xt}+\frac{2}{t}u_{t}=(u_{t}+u_{x})(u_{t}+t^{4}u_{t}^{3}) \end{equation}
		
		with the initial value conditions $u(x,1)=x$ and $u_{t}(x,1)=\frac{1}{\sqrt{2e^{-2x}-1}}$.
		
	\end{Example}
	
	\begin{remark} In ~\citep{Mahouton}, the author studied special case of (\ref{eq25}) where $u$ is a single variable function of $t$.
		\\
		In (\ref{eq25}), if we suppose that $u$ is only a function of $t$, we get 
		$$u''(t)+\frac{2}{t}u'=u'^{2}+(tu')^{4}$$ which is the differential equation (61) investigated by the authors in ~\citep{Mahouton}.
	\end{remark}
	
	\textbf{Solution:}	
	
	(\ref{eq25}) is reduced to the differential equation 
	\begin{equation}\label{eq26}u_{t}(x,t)=K(u)\frac{1}{t^{2}}
	\end{equation}
	where $K$ is the general solution of $$K'(u)=k(u)+u^{3}.$$
	$$K(u)=\pm(A e^{-2u}-1)^{-\frac{1}{2}}$$
	Then the general solutions of \ref{eq25} are given by 
	$$u_{t}(x,t)=t^{-2}(f(x-t)e^{-2u}-1)^{-\frac{1}{2}}$$ where
	$f$ is an arbitrary function.
	
	By generating Mathematica codes, we obtain implicit solutions of (\ref{eq25}):
	
$$	\sqrt{2 {\mathrm e}^{-2 u\left(x,t\right)}-1} \,  t-\arctan\! \left(\sqrt{2 {\mathrm e}^{-2 u\left(x,t\right)}-1}\right) t-t \sqrt{2 {\mathrm e}^{-2 x}-1}+t \arctan\! \left(\sqrt{2 {\mathrm e}^{-2 x}-1}\right)+t-1=0$$
	
\begin{figure}[h]
		\centering
		\begin{subfigure}[b]{0.5\textwidth}
			\centering
			\includegraphics[width=\textwidth]{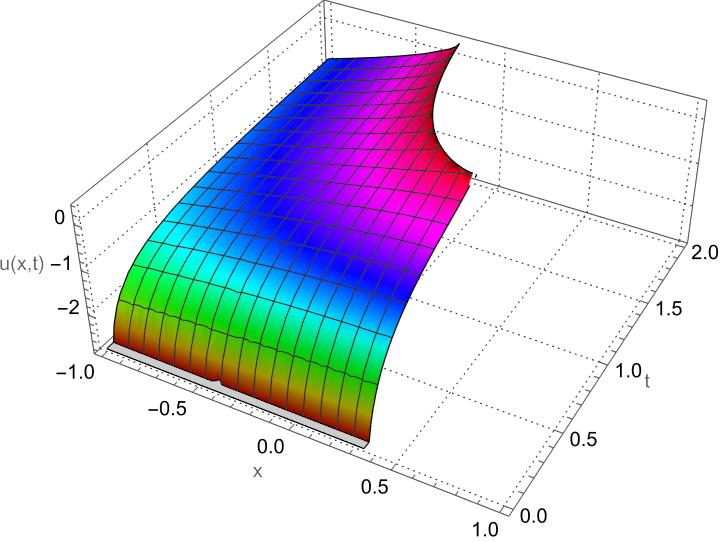}
			\caption{ }
		\end{subfigure}
		\hfill
		\begin{subfigure}[b]{0.4\textwidth}
			\centering
			\includegraphics[width=\textwidth]{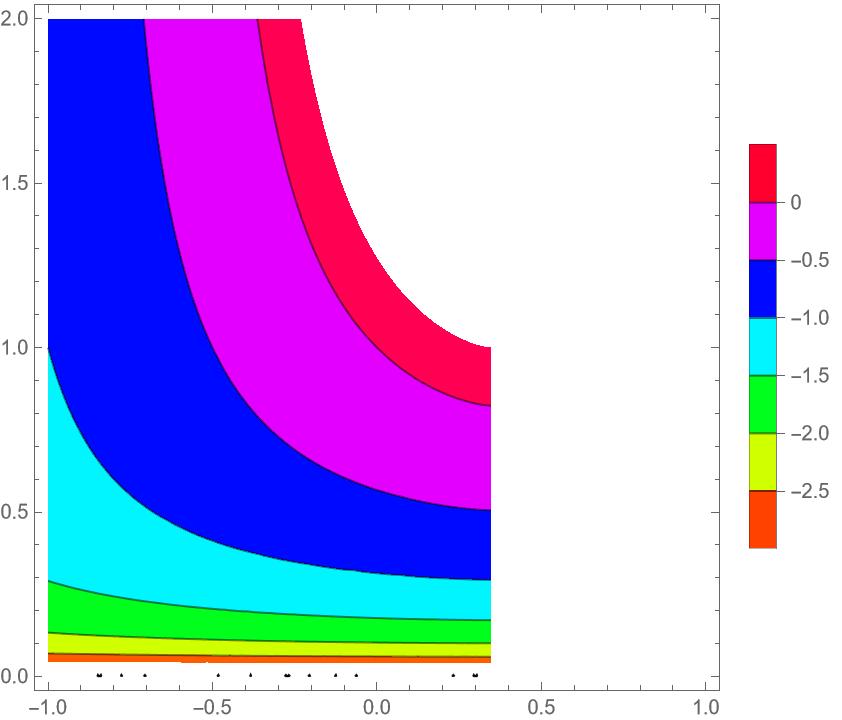}
			\caption{ }
		\end{subfigure}
		\hfill
		\begin{subfigure}[b]{0.5\textwidth}
			\centering
			\includegraphics[width=\textwidth]{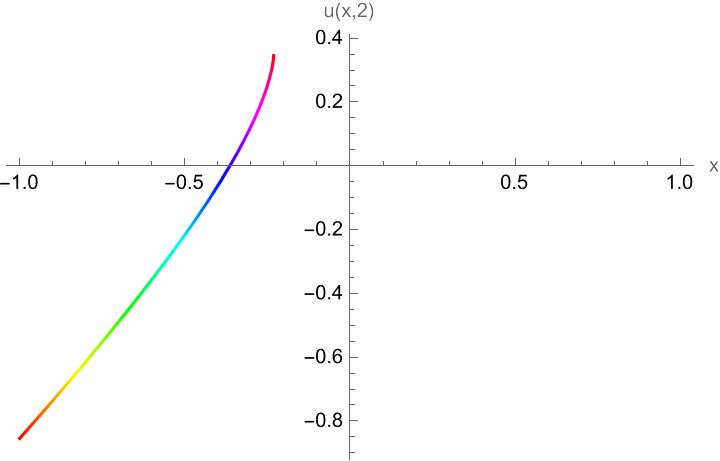}
			\caption{}
		\end{subfigure}
		\caption{The profile of the solutions in  (\ref{eq25}) with $u(x,1)=x$ and $u_{t}(x,1)=\frac{1}{\sqrt{2e^{-2x}-1}}$: (a) and (b) 3D and Contour plots  with $-1\leq x\leq1$ and $0\leq t\leq 2$, (c) 2D plot at $t=2$.}
		\label{fig7}
	\end{figure}

	Figure \ref{fig7} shows the 2D, 3D and contour plots of the solutions in (\ref{eq25}) within $-1\leq x\leq1$ and $0\leq t\leq 2$ for 3D and contour graphs, $t=2$ for 2D graph. 
	
	\begin{Example} Let $m=0$, $b(t)=-\frac{1}{t}$.
		Suppose that function $F$ satisfies $F(s,w)= 1+2\frac{s}{w}$.
		\begin{equation}\label{eq27} u_{tt}+au_{xt}-\frac{1}{t}u_{t}=t(u_{t}+au_{x})(1+2t\frac{u}{u_{t}}). \end{equation}
	\end{Example}
	\begin{remark} 
		A special case of our findings was recorded in ~\citep{Mahouton}.
		If we suppose that $u$ is only a function of $t$ in (\ref{eq27}) , we get 
		$$u''(t)-(\frac{1}{t}+t)u'-2t^{2}u=0,$$ which is exactly the differential equation (78) investigated by the authors of ~\citep{Mahouton}.
	\end{remark}
	
	\textbf{Solution:}	
	
	(\ref{eq27}) is reduced to the differential equation 
	\begin{equation}\label{eq28}u_{t}(x,t)=tK(u).
	\end{equation}
	where $K$ satisfies \begin{equation}\label{eq29}K'(u)K(u)=K(u)+2u.\end{equation}
	Two particular solutions to (\ref{eq29}) are given by $K_{1}(u)=2u$ and $K_{2}(u)=-u$.\\
	The general solutions of (\ref{eq29}) satisfy the algebraic equation  
	\begin{equation}\label{eq30}(K(u)-2u)^{2}(K(u)+u)=A.\end{equation} where $A$ is an arbitrary constant.\\
	A real solution of  (\ref{eq30}) can be computed to yield   
	$$K(u)=\frac{\sqrt[3]{\sqrt{A^2-4 A u^3}-2 u^3+A}}{\sqrt[3]{2}}+\frac{\sqrt[3]{2} u^2}{\sqrt[3]{\sqrt{A^2-4 A u^3}-2 u^3+A}}+u$$
	If $a=1$, $K(u)=\phi(u,x-t)$ where $\phi$ is an arbitrary function.
	
	Hence, we obtain an implicit solution of (\ref{eq27}) as $$u_{t}(x,t)=t\phi(u,x-t)$$ 
	
	We checked the implicit solutions of (\ref{eq27}) by generating Mathematica codes considering the initial conditions  $u(x,1)=1$. 
	
\begin{figure}[h]
		\centering
		\begin{subfigure}[b]{0.5\textwidth}
			\centering
			\includegraphics[width=\textwidth]{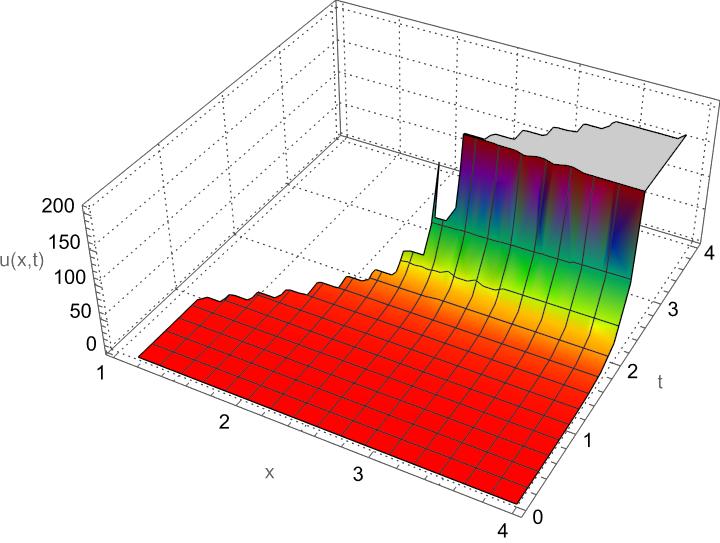}
			\caption{ }
		\end{subfigure}
		\hfill
		\begin{subfigure}[b]{0.4\textwidth}
			\centering
			\includegraphics[width=\textwidth]{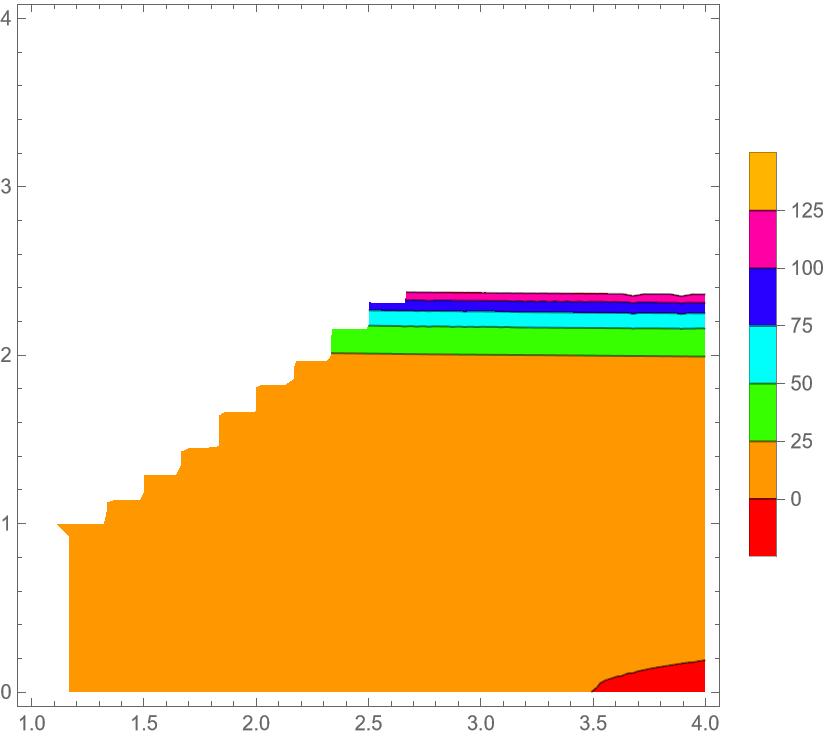}
			\caption{ }
		\end{subfigure}
		\hfill
		\begin{subfigure}[b]{0.5\textwidth}
			\centering
			\includegraphics[width=\textwidth]{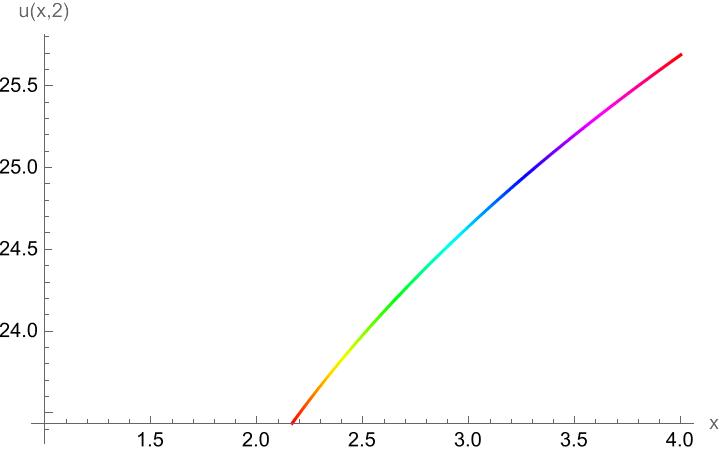}
			\caption{}
		\end{subfigure}
		\caption{The profile of the solutions in  (\ref{eq27}) with $u(x,1)=1$ and $u_{t}(x,1)=x^{2}$: (a) and (b) 3D and Contour plots  with $1\leq x\leq 4$ and $0\leq t\leq 4$, (c) 2D plot at $t=2$.}
		\label{fig8}
	\end{figure}

	Figure \ref{fig8} shows the 2D, 3D and contour plots of the solutions in (\ref{eq27}) within $1\leq x\leq4$ and $0\leq t\leq 4$ for 3D and contour graphs, $t=2$ for 2D graph. 
	
	\section{Conclusion}
	
	In this paper, we presented a new method, a combination of the variation of parameters and other techniques, such as the method of characteristics, to derive exact solutions of nonlinear partial differential equations alongside specific initial conditions, a framework extensively applied in mathematical physics. Illustrative examples were provided to demonstrate the applicability of this method. Problems that are nontrivial when approached with conventional methods now appear straightforward, as the resulting functions are univariate. Our research findings indicate that fusing established classical techniques with innovative approaches yields efficient analytical solutions.
	
	\bibliographystyle{unsrtnat}
	\bibliography{references}  

\end{document}